\documentclass{article}
\usepackage{cmap}
\usepackage[T2C]{fontenc}
\usepackage[utf8]{inputenc}
\usepackage[english,russian]{babel}
\usepackage[tbtags]{amsmath}
\usepackage{amsfonts,amssymb}
\usepackage[matrix,arrow,curve]{xy}
\usepackage{mathrsfs}
\usepackage[fanpaper]{faa_authors-utf8_eng}
\usepackage{hypbmsec}

\usepackage{tabularx}
\usepackage{multirow}
\usepackage{rotating}
\usepackage{hhline}
\usepackage{esint}

\newenvironment{bad}{\spaceskip=0.33emplus0.6emminus0.15em
\immediate\write5{\string\bad}}{}

\theoremstyle{plain}

\newtheorem{theorem}{Theorem}
\newtheorem{proposition}{Proposition}
\newtheorem{lemma}{Lemma}

\theoremstyle{definition}
\newtheorem{proof}{Proof}
\newtheorem{definition}{Definition}
\newtheorem{remark}{Remark}

\newtheorem{example}{Example}

\begin{document}
\selectlanguage{english}
\def\refname{Literature}

\rusisspages{15--39}

\frenchspacing

\title{Rings of invariants of three-dimensional upper-triangular representations of finite groups}

\author{ABDULKADYR BUCHAEV}

\address{Moscow Institute of Physics and Technology, Dolgoprudny, Russia}
\email{buchaev.aia@phystech.edu}



\markboth{Invariants of upper triangular groups}{A.\,Y.~Buchaev}
\maketitle

\begin{fulltext}

\thanks{The research is supported by the MSHE ``Priority 2030'' strategic academic leadership program.}


\begin{abstract}
The work proves that, for finite three-dimensional upper triangular groups over a field of odd characteristic generated by pseudoreflections, maximal unipotent subgroup of which is abelian, the ring of invariants is polynomial if and only if the maximal unipotent subgroup is generated by pseudoreflections or does not contain transvections.
\end{abstract}


\section{Introduction}
\label{s1}

One of the problems addressed in algebraic geometry is the computation of quotient varieties under group actions. Since this problem is very complex in its general form, it is interesting to begin by solving it for groups or varieties of a special type. This paper deals with the action of a finite linear group $G \subset \operatorname{GL}(W)$, where $W$ is a vector space of dimension $n$ over a field $\mathbb{K}$. Here and thereafter, the field $\mathbb{K}$ is assumed to be algebraically closed.

The left action of the group $G$ on $W$ induces a left action of $G$ on the symmetric algebra of the dual space, which is isomorphic to the ring of regular functions $\mathbb{K}[x_1, \ldots, x_n]$; this action is defined as follows:
\begin{gather}
 f(x) \in \mathbb{K}[W], 
 \qquad g\in G,
 \label{eq1}
\\
 g\colon  f(x) \mapsto f(g^{-1} \cdot x).
 \label{eq2}
\end{gather}

The ring of invariants under the action of the group $G$ is then defined:
\begin{equation}
 \label{eq3}
 \mathbb{K}[W]^G=\{f(x) \in \mathbb{K}[W]\colon  g\cdot f=f\}.
\end{equation}
It is easy to see that the ring of invariants is a graded subalgebra of $\mathbb{K}[W]$, which reduces the problem of finding invariants to the problem of finding homogeneous invariants. As is known, if $G$ is a finite group, then $\mathbb{K}[W]^G$ is a finitely generated $\mathbb{K}$-algebra.

We will be interested in the question of the non-singularity of the quotient variety. It is well known that in the case where $W$ is a finite-dimensional vector space and $G$ is a finite linear group, the non-singularity of the quotient variety $W/G$ is equivalent to the polynomial property of the ring $\mathbb{K}[W]^G$ (that is, $\mathbb{K}[W]^G$ is a free $\mathbb{K}$-algebra).

\begin{lemma}[{\rm (see\ \cite[Lemma 1]{1})}]
\label{l1} Suppose a finite group $G$ acts linearly on a vector space $W$ of finite dimension over a field $\mathbb{K}$. Then the ring of invariants $\mathbb{K}[W]^G$ is polynomial if and only if the variety $W/G$ is non-singular at the point $0$.
\end{lemma}

The starting point of the research is the theorem of  Shephard and Todd     \cite{2}, Chevalley~\cite{3}, and Serre~\cite{1}. It is known as the Chevalley--Shephard--Todd theorem. It provides an answer to the question of the polynomiality of the ring of invariants in the case where the characteristic of the field $\mathbb{K}$ does not divide the order of the group $G$.

\begin{definition}
\label{d1}
Let $W$ be a vector space of dimension $n$ over a field $\mathbb{K}$. A linear transformation $A \in \operatorname{GL}(W)$ of finite order is called a \textit{pseudoreflection} if
\begin{equation}
\label{eq4}
\dim \ker (A-E)=n-1,
\end{equation}
where $E$ is the identity transformation.
\end{definition}

In other words, $A \in \operatorname{GL}(W)$ is a pseudoreflection if $A \not= E$ and there exists a hyperplane $H \subset W$ such that $A$ acts identically on $H$, and $A$ has finite order.
Also, for convenience, we will sometimes consider the identity operator to be a pseudoreflection.
If a pseudoreflection $A$ is diagonalizable, then $A$ is called a \textit{homology}. Otherwise, $A$ is called a \textit{transvection}.

\begin{theorem}[{\rm (Chevalley--Shephard--Todd; see\ \cite{1}--\cite{3})}]
\label{t1} Let $W$ be a vector space of finite dimension over a field $\mathbb{K}$ of characteristic $p$, $p \ge 0$, and let $G \subseteq \operatorname{GL}(W)$ be a finite subgroup. Then, if the ring of invariants $\mathbb{K}[W]^G$ is polynomial, the group $G$ is generated by pseudoreflections. If the field characteristic $p$ does not divide the order of the group $G$, then the converse is also true.
\end{theorem}

However, if the field characteristic $p$ divides $|G|$, i.e., in the so-called modular case, the converse statement is false -- there exist groups generated by pseudoreflections whose ring of invariants is not polynomial (see\ \cite[p.~8.2]{4} and \cite[Proposition 3.1]{5}). 
Therefore, in the modular case, the question remains open: what could be the criterion for the polynomiality of the ring of invariants?

Kemper and Malle answered this question in the case where $G$ acts on $W$ via an irreducible representation.

\begin{theorem}[{\rm (Kemper--Malle; see\ \cite{5})}]
\label{t2}
Let $W$ be a finite-dimensional vector space, and let $G \subseteq \operatorname{GL}(W)$ be a finite irreducible group. Then $\mathbb{K}[W]^G$ is polynomial if and only if $G$ is generated by pseudoreflections and the pointwise stabilizer $G_U \subseteq G$ of any nontrivial subspace $U$ in $W$ has a polynomial ring of invariants $\mathbb{K}[W]^{G_U}$.
\end{theorem}

\begin{remark}
\label{r1}
In their original article \cite{5}, the authors made an error by formulating the theorem for the symmetric algebra $S(W)$ instead of the ring of regular functions $\mathbb{K}[W]$, which is isomorphic to $S(W^*)$. Later, in the book by Kemper and Derksen \cite[Theorem 3.9.8]{6}, the inaccuracy was corrected.
\end{remark}

Also, in 2023, Braun announced a generalization of the result of Kemper and Malle to subgroups of the group $\operatorname{SL}(W)$, which are not necessarily irreducible:

\begin{theorem}[{\rm (Braun; see\ \cite{7})}]
\label{t3}
Let $W$ be a finite-dimensional vector space, and let $G  \subseteq\operatorname{SL}(W)$ be a finite group. Then $\mathbb{K}[W]^G$ is polynomial if and only if $G$ is generated by pseudoreflections and the pointwise stabilizer $G_U$ in $G$ of any nontrivial subspace $U$ in $W$ has a polynomial ring of invariants $\mathbb{K}[W]^{G_U}$.
\end{theorem}

Nevertheless, all known counterexamples to the converse of Theorem \ref{t1} in the modular case that do not fall under the Kemper-Malle criterion (i.e., reducible indecomposable representations) start from dimension 4, and in dimension 2 the theorem holds in both directions. This raised the question: is there a ring of invariants of a reducible indecomposable representation in dimension 3 that is non-polynomial? The following theorem, proven in the present work, answers this question:

\begin{theorem}
\label{t4}
Let $G$ be a finite group acting on a three-dimensional vector space $W$ over a field of odd characteristic $p$ by upper triangular matrices, and $G$ is generated by pseudoreflections. Let $H$ be the subgroup of $G$ consisting of all elements acting by unipotent matrices. If $H$ is abelian, then $\mathbb{K}[W]^G$ is polynomial if and only if $H$ is generated by pseudoreflections or does not contain transvections.
\end{theorem}

Indeed, below we construct an example of a group $G$ satisfying the hypotheses of the theorem such that its maximal unipotent subgroup $H$ is abelian and contains transvections, but is not generated by them. The theorem then implies that $\mathbb{K}[W]^G$ is not polynomial. Moreover, $G$ acts on $W$ via a reducible indecomposable representation, thus providing the desired example.

We shall prove the theorem after first establishing a series of technical results.

\section{Preliminaries}
\label{s2}

In this section, we will describe the main lemmas that will be useful in studying the question of the polynomiality of the ring of invariants. We use standard facts from the theory of invariants of finite groups; see, for example, \cite{4} and \cite{8}.

\begin{theorem}[{\rm (see\ \cite[Theorem 1.3.1]{8})}] 
\label{t5}
Let $\mathbb{K}$ be a field, $A$ a finitely generated commutative $\mathbb{K}$-algebra, and $G$ a finite group acting on $A$ by automorphisms. Then $A^G$ is a finitely generated algebra, and $A$ is a finitely generated module over $A^G$ (in particular, $A$ is integral over $A^G$).
\end{theorem}

\begin{proposition}
\label{p1}
Let a finite group $G$ act linearly on a vector space $W$ of dimension $n$ over a field $\mathbb{K}$, and suppose $\mathbb{K}[W]^G$ is a polynomial $\mathbb{K}$-algebra. Then $\mathbb{K}[W]^G$ has exactly $n$ algebraically independent generators.
\end{proposition}

We state two lemmas that form the basis for verifying whether a set of $n$ homogeneous polynomials is an algebraically independent system of generators for the algebra $\mathbb{K}[W]^G$.

\begin{definition}
\label{d2}
Let $R$ be a finitely generated graded $\mathbb{K}$-algebra. Then a set of homogeneous elements $\{f_1, \ldots, f_n\}$ from $R$, where $n$ is the Krull dimension of $R$, such that $R$ is finitely generated as a module over $\mathbb{K}[f_1, \ldots, f_n]$, is a homogeneous system of parameters for $R$.
\end{definition}

Note that the elements of a homogeneous system of parameters are algebraically independent over $\mathbb{K}$. Next, we will state a lemma that allows us to determine whether a given set of homogeneous elements is a homogeneous system of parameters.

\begin{lemma}[{\rm (see\ \cite[Lemma 2.6.3]{4})}]
\label{l2}
Let $A \subset \mathbb{K}[W]$ be a Noetherian graded subring of dimension $n = \dim W$. Suppose $\{h_1, \ldots, h_n\}$ is a system of homogeneous elements of the ring $A$. Then $\{h_1, \ldots, h_n\}$ is a homogeneous system of parameters for $A$ if and only if the common zero set of the polynomials $h_1, \ldots, h_n$ is $\{0\}$.
\end{lemma}

Finally, we state a lemma that allows us to determine whether the elements of a homogeneous system of parameters generate $\mathbb{K}[W]^G$.

\begin{lemma}[{\rm (see\ \cite[Corollary 3.1.6]{4})}]
\label{l3}
Let $\{f_1, \ldots, f_n\}$ be a homogeneous system of parameters for $\mathbb{K}[W]^G$ of degrees $d_1, \ldots, d_n$, respectively. Then
\begin{equation}
\label{eq5}
\prod_{i=1}^n d_i=|G| \Longleftrightarrow \mathbb{K}[W]^G=\mathbb{K}[f_1, \dots, f_n].
\end{equation}
\end{lemma}

Thus, $\mathbb{K}[W]^G$ is a polynomial algebra if and only if there exist homogeneous elements
\begin{equation}
\label{eq6}
 f_1, \dots, f_n \in \mathbb{K}[W]^G,
 \qquad  n=\dim W,
\end{equation}
such that their common zero set is $\{0\}$ and the product of their degrees equals the order of the group $G$.

Indeed, if such $f_1, \ldots, f_n$ exist, then by Lemma~\ref{l2} they form a homogeneous system of parameters and are algebraically independent over $\mathbb{K}$. Applying Lemma~\ref{l3}, we then conclude that they generate $\mathbb{K}[W]^G$. Conversely, suppose that $\mathbb{K}[W]^G$ is a polynomial $\mathbb{K}$-algebra. Choose a set of generators $f_1, \ldots, f_n$ of $\mathbb{K}[W]^G$; by Proposition~\ref{p1}, there are exactly $n = \dim W$ of them. By definition, $f_1, \ldots, f_n$ form a homogeneous system of parameters for $\mathbb{K}[W]^G$. Hence, by Lemma~\ref{l2}, their common zero set is $\{0\}$. Moreover, since $f_1, \ldots, f_n$ generate $\mathbb{K}[W]^G$ and form a homogeneous system of parameters, Lemma~\ref{l3} implies that the product of their degrees is equal to $|G|$.

\section{Linearization of quotient group action}
\label{s3}

Often, when studying the ring of invariants of a group $G$, it is useful to select a normal subgroup $H$ within $G$, study its ring of invariants, and then, knowing the ring of invariants $\mathbb{K}[W]^H$, recover certain properties of $\mathbb{K}[W]^G$. In this section, we will prove a statement that will allow us to carry out the described transition:

\begin{proposition}
\label{p2} 
{\bad
Let $G$ be a finite linear group acting on a vector space $W$ of dimension $n$ over a field $\mathbb{K}$ of characteristic $p >\nobreak 0$,} and let $G$ be generated by pseudoreflections. Let $H$ be its normal subgroup such that the exponent of $p$ in $|H|$ is equal to the exponent of $p$ in $|G|$. Then, if $\mathbb{K}[W]^H$ is polynomial, $\mathbb{K}[W]^G$ is polynomial.
\end{proposition}

This statement will allow us to reduce questions about the polynomiality of rings of invariants under group actions to questions about the polynomiality of rings of invariants of their Sylow $p$-subgroups, where $p$ is the characteristic of the base field.

Obviously, if $H$ is a normal subgroup of $G$, then
 \begin{equation}
 \label{eq7}
 \mathbb{K}[W]^G=(\mathbb{K}[W]^H)^{G/H}=\mathbb{K}[W/H]^{G/H}.
 \end{equation}

The idea is that if $|H|$ has the same exponent of $p$ as $|G|$, then we can apply the Chevalley-Shephard-Todd theorem to the non-modular group $G/H$ acting on the space $W/H$, which is isomorphic to an affine space of dimension $n$ since $\mathbb{K}[W]^H$ is polynomial. However, the problem is that the induced action $\rho$ of $G/H$ on $W/H$ may be nonlinear. To establish the polynomiality of the ring of invariants, we will consider another action $\rho'$ of $G/H$, which is linear, and the group $G/H$ under this action is generated by pseudoreflections; moreover, locally this action coincides with the action $\rho$, from which the required result will follow.

First, note that by Lemma~\ref{l1}, it suffices to prove that the ring $\mathbb{K}[W]^G$ is regular at the point 0. Moreover, $\mathbb{K}[W]^G$ is regular at zero if and only if the completion $\widehat{\mathbb{K}[W]^G_\mathfrak{m}}$ of the local ring $\mathbb{K}[W]^G_\mathfrak{m}$ is regular, where $\mathfrak{m}$ is the maximal ideal corresponding to the point 0.
Also
\begin{equation}
\label{eq8}
	\widehat{\mathbb{K}[W]^G_\mathfrak{m}}=\widehat{\mathbb{K}[W/H]^{G/H}_\mathfrak{m}},
\end{equation}
so it suffices to prove that the ring $\widehat{\mathbb{K}[W/H]^{G/H}_\mathfrak{m}}$ is regular.

For this, we will need the following lemma:

\begin{lemma}[{\rm (see\  \cite[Lemma 2.3]{9})}]
\label{l4}
Let $R = \mathbb{K}[[x_1, x_2, \ldots, x_n]]$ be the local ring of formal power series in $n$ variables over a field $\mathbb{K}$, and let $G$ be a finite group acting on $R$ by local automorphisms. Suppose the characteristic of the field $\mathbb{K}$ does not divide the order of the group $G$. Then one can choose local parameters $y_1, \ldots, y_n$ in $R$ such that $G$ acts by linear substitutions on $y_1, \ldots, y_n$.
\end{lemma}

Note that by assumption, $W/H$ is isomorphic to an affine space, so the completion of the ring $\mathbb{K}[W/H]$ at zero is isomorphic to the ring of formal power series in $n = \dim W/H$ variables:
 \begin{equation}
 \label{eq9}
 \widehat{\mathbb{K}[W/H]_\mathfrak{m}} \simeq \mathbb{K}[[x_1, x_2, \dots, x_n]],
 \qquad n=\dim W/H.
 \end{equation}
    Also, by assumption, $|G/H|$ is not divisible by the characteristic of the field $\mathbb{K}$, so Lemma~\ref{l4} is applicable to the ring $\widehat{\mathbb{K}[W/H]_\mathfrak{m}}$.
    
    The group $G/H$ acts nonlinearly on the space $W/H$, but according to Lemma~\ref{l4}, in $\widehat{\mathbb{K}[W/H]_\mathfrak{m}}$ one can choose a system of generators $y_1, \ldots, y_n$ such that the action of $G/H$ on them is linear. Inducing the action of $G/H$ on the linear span $\langle y_1, \dots, y_n \rangle$ of series $y_1, \ldots, y_n$ over $\mathbb{K}$, we obtain an induced action of $G/H$ on the dual space $V=\langle y_1, \dots, y_n \rangle^*$ and $\mathbb{K}[V]=\mathbb{K}[y_1, \dots, y_n], \widehat{\mathbb{K}[V]_\mathfrak{m}}= \mathbb{K}[[y_1, \dots, y_n]]$.

    Then the actions of $G/H$ on $\widehat{\mathbb{K}[W/H]_\mathfrak{m}}$ and $\widehat{\mathbb{K}[V]_\mathfrak{m}}$ are equivalent (i.e., $\widehat{\mathbb{K}[W/H]_\mathfrak{m}}$ and $\widehat{\mathbb{K}[V]_\mathfrak{m}}$ are isomorphic as $\mathbb{K}[G/H]$-algebras). Consequently, to prove that the ring $\widehat{\mathbb{K}[W/H]^{G/H}_\mathfrak{m}}$ is regular at the point 0, it suffices to show that $\widehat{\mathbb{K}[V]_\mathfrak{m}^{G/H}}$ is regular at the point 0, which, in turn, is equivalent to the ring $\mathbb{K}[V]^{G/H}$ being polynomial. To do this, we will use Chevalley--Shephard--Todd theorem: the action of $G/H$ on $V$ is nonmodular, but we need to show that $G/H$ is generated by pseudoreflections.

    To achieve this, we need to formulate what a pseudoreflection is in terms of the ring of regular functions of the space. This is provided by the following lemma:

\begin{lemma}
\label{l5}
Let $V$ be a vector space of dimension $n$ over an algebraically closed field $\mathbb{K}$, $\mathbb{K}[V] = \mathbb{K}[x_1,\ldots, x_n]$, and let $g \in\operatorname{GL}(V)$ be a linear operator of finite order. Then $g$ is a pseudoreflection if and only if there exists an irreducible $\alpha \in \mathfrak{m} = (x_1, \ldots, x_n)$ such that for any polynomial $f \in \mathbb{K}[V]$, we have $g \cdot f - f \in (\alpha)$.
\end{lemma}

\begin{proof} Indeed, if $g$ is a pseudoreflection and $\alpha$ is a linear polynomial defining the hyperplane on which $g$ acts identically, then for each point $x$ of this hyperplane, $g \cdot x = x$, whence $g \cdot f - f$ is identically zero on this hyperplane. Therefore, by Hilbert's Nullstellensatz, there exists $k$ such that $(g \cdot f - f)^k \in (\alpha)$, but from the irreducibility of $\alpha$ it follows that $g \cdot f - f \in (\alpha)$.

Conversely, if for some irreducible polynomial $\alpha \in \mathfrak{m}$ it is true that for any $f \in \mathbb{K}[V]$ the polynomial $g \cdot f - f$ lies in $(\alpha)$, then taking $f$ to be the coordinate functions $x_i$, we find that $g \cdot f$ coincides with $f$ on the set $V(\alpha)$ of zeros of the polynomial $\alpha$. Consequently, $g$ does not change any of the coordinates of points in this set, meaning $g$ acts identically on it. But by Krull's principal ideal theorem, the nonempty closed set $V(\alpha)$ has codimension 1, and the set of fixed points of $g \in GL(V)$ is a vector subspace of $V$, whence we conclude that the set of fixed points of the operator $g$ contains a hyperplane. Therefore, $g$ is a pseudoreflection. \hfill $\square$
\end{proof}

\begin{remark}
\label{r2}
Note that for the converse statement of the lemma, it suffices to require that $g \cdot f - f$ lie in $(\alpha)$ for the coordinate functions $x_1, \ldots, x_n$.
\end{remark}

Note that in the previous lemma, the polynomial $\alpha$ need not be linear -- this allows generalizing the concept of a pseudoreflection to the case of a nonlinear action.
    
\begin{definition}
\label{d3}
In the framework of this work, we will call an automorphism $g$ of an affine space $V$ of dimension $n$ of finite order a pseudoreflection if the set of fixed points of the automorphism $g$ contains the point $0$ and has codimension 1 or, equivalently, if there exists an irreducible polynomial $\alpha \in \mathfrak{m} = (x_1, \ldots, x_n)$ such that for any $f \in \mathbb{K}[V]$, we have $g \cdot f - f \in (\alpha)$.
\end{definition}

Next, we will show that if $G$, acting on $W$, is generated by pseudoreflections $g_1, \ldots, g_k$, then their images $g_1H, \ldots, g_kH$ in the quotient group $G/H$ acting on the space $W/H \simeq \mathbb{A}^n$ are also pseudoreflections (and, of course, generate $G/H$).

\begin{lemma}
\label{l6}
Let $W$ be a vector space of dimension $n$ over an algebraically closed field $\mathbb{K}$, let the group $G  \subseteq \operatorname{GL}(W)$ be generated by pseudoreflections, and let $H$ be a normal subgroup of $G$ such that $W/H \simeq \mathbb{A}^n$. Then the (nonlinear) quotient action of $G/H$ on $W/H$ is also generated by pseudoreflections.
\end{lemma}

\begin{proof} First, note that the point $0 \in W$ is fixed by the entire group $G$, whence we conclude that $0 \in W/H$ is fixed by the entire $G/H$. Therefore, it suffices to prove that for any $i$, the codimension of the set of fixed points of the automorphism $g_iH$ is at most 1.
    
    Next, observe that the quotient morphism $W \twoheadrightarrow W/H$ is a finite morphism: indeed, it corresponds to the embedding of rings $\mathbb{K}[W]^H \hookrightarrow \mathbb{K}[W]$, and $\mathbb{K}[W]$ is integral over $\mathbb{K}[W]^H$, so the quotient morphism is finite. Consequently, the quotient $W \twoheadrightarrow W/H$ is a closed morphism, which, moreover, preserves the dimensions of closed subvarieties.

    Let $L$ be the set of fixed points of the automorphism $g \in G$. Then, since the group $H$ is normal, for $x \in L$ we have
 \begin{equation}
 \label{eq10}
 gHx=Hgx=Hx,
 \end{equation}
whence it follows that $Hx$ is a fixed point of the automorphism $gH \in G/H$ acting on $W/H$. Therefore, $L$ maps under the quotient morphism to the set of fixed points of the automorphism $gH$, and
$$
 \operatorname{codim}_W L=\operatorname{codim}_{W/H} L/H.
$$
Thus, if $g$ is a pseudoreflection on $W$, then $gH$ is a pseudoreflection on $W/H$. \hfill $\square$
\end{proof}
    
Therefore, if $G$ is generated by pseudoreflections on $W$, then $G/H$ is generated by pseudoreflections on $W/H$. It remains only to show that the same holds after the change of basis in Lemma~\ref{l4}.

\begin{lemma}
\label{l7}
Let $V$ and $V'$ be affine spaces of dimension $n$ over an algebraically closed field $\mathbb{K}$, with $\mathbb{K}[V] = \mathbb{K}[x_1, \ldots, x_n]$, $\mathbb{K}[V'] = \mathbb{K}[y_1, \ldots, y_n]$, and let $g \in \operatorname{Aut}(V)$, $g' \in \operatorname{GL}(V')$. Then, if the induced actions of $g$ and $g'$ on formal power series rings $\mathbb{K}[[V]]= \mathbb{K}[[x_1, \dots, x_n]]$ and $\mathbb{K}[[V']]=\mathbb{K}[[y_1, \dots, y_n]]$, respectively, are equivalent and $g$ is a pseudoreflection, then $g'$ is also a pseudoreflection.
\end{lemma}

\begin{proof} Since $g$ is a pseudoreflection on $V$, by Lemma \ref{l5} there exists a polynomial $\alpha \in \mathfrak{m} = (x_1, \ldots, x_n)$ such that for any $f \in \mathbb{K}[V]$, we have $g \cdot f - f \in (\alpha)$. Note that upon localization at the maximal ideal $\mathfrak{m}$, this inclusion still holds. We will show that the same condition holds in the ring of formal power series -- the completion of the localization at the maximal ideal.

Indeed, let $f$ be an arbitrary series and $f_k$ be the sum of its monomials of degree at most $k$. Then $g \cdot f_k - f_k \in (\alpha)$ for all integers $k$. Moreover, $(f-f_k) \in \mathfrak{m}^{k+1}$, and hence $g\cdot(f-f_k) \in \mathfrak{m}^{k+1}$ since $g$ does not lower the grading; consequently, $g\cdot (f-f_k)-(f-f_k) \in \mathfrak{m}^{k+1}$. Then, upon passing to the quotient by the power of the maximal ideal, we obtain that $g \cdot f - f \in ((\alpha)+ \mathfrak{m}^{k+1})/\mathfrak{m}^{k+1}$ for any positive integer $k$, and hence $g \cdot f - f \in (\alpha)$.

Let $\varphi \colon  \mathbb{K}[[x_1, \dots, x_n]] \simeq \mathbb{K}[[y_1, \dots, y_n]]$ be an isomorphism such that the actions of $g$ and $g'$ are equivalent, i.e., $\varphi g=g' \varphi$. Set $\alpha' = \varphi(\alpha)$. Note that $\alpha'$ is a formal power series, which is not necessarily a polynomial. Observe that
 \begin{equation}
 \label{eq11}
 \alpha \in (x_1, \dots, x_n)
 \quad \Longleftrightarrow\quad \alpha' \in (y_1, \dots, y_n).
 \end{equation}

        Then we have that
 \begin{equation}
 \label{eq12}
    \forall\, f \in \mathbb{K}[[V]]\quad g \cdot f -f \in (\alpha) \quad \Longleftrightarrow\quad \forall\, f \in \mathbb{K}[[V]]\quad g' \cdot \varphi(f) - \varphi(f) \in (\alpha'),
 \end{equation}
        but since $\varphi$ is an isomorphism, it follows that
 \begin{equation}
 \label{eq13}
 \forall\, f' \in \mathbb{K}[[V']] \quad g' \cdot f' - f' \in (\alpha').
 \end{equation}

        In particular, for all $i$ from 1 to $n$, we have
 \begin{equation}
 \label{eq14}
 g' \cdot y_i - y_i \in (\alpha').
 \end{equation}
        Since $g' \in \operatorname{GL}(V')$, polynomial $g' \cdot y_i - y_i$ is linear, whence
\begin{equation}
 \label{eq15}
g' \cdot y_i - y_i=c_i \alpha'_1,
\end{equation}
        where $c_i \in \mathbb{K}$, and $\alpha'_1$ is degree-1 form of $\alpha'$, $\alpha'_1 \in \mathbb{K}[V']$.
        From this and Remark~\ref{r2}, it follows that $g'$ is a pseudoreflection on $V'$. \hfill $\square$
\end{proof}

    Finally, we combine all the statements in the proof of Proposition~\ref{p2}.

\begin{proof}[of Proposition~\ref{p2}]
         The ring $\mathbb{K}[W]^G$ is polynomial if and only if its localization $\mathbb{K}[W]^G_\mathfrak{m}$ at the maximal ideal $\mathfrak{m}$ corresponding to the point 0 is regular. This is equivalent to the completion $\widehat{\mathbb{K}[W]^G_\mathfrak{m}}$ of this localization being regular. On the other hand, from the equality $\mathbb{K}[W]^G = \mathbb{K}[W/H]^{G/H}$, we obtain that
 \begin{equation}
  \label{eq16}
     \widehat{\mathbb{K}[W/H]^{G/H}_\mathfrak{m}} \simeq \widehat{\mathbb{K}[W]^G_\mathfrak{m}},
 \end{equation}
 therefore, it suffices to show that $\widehat{\mathbb{K}[W/H]^{G/H}_\mathfrak{m}}$ is regular.
 To this end, consider the linear action of $G/H$ on a vector space $V$, $V^* = \langle y_1, \ldots, y_n \rangle$, where $y_1, \ldots, y_n$ are the local parameters provided by Lemma~\ref{l4} for $R=\widehat{\mathbb{K}[W/H]_\mathfrak{m}}$ and the group $G/H$. Observe that $\widehat{\mathbb{K}[V]_\mathfrak{m}}$ and $\widehat{\mathbb{K}[W/H]_\mathfrak{m}}$ are isomorphic as $(G/H)$-modules. By Lemmas~\ref{l6} and~\ref{l7}, the group $G/H$ is generated by pseudoreflections on $V$. Moreover, the condition $\operatorname{ord}_p(|H|)=\operatorname{ord}_p(|G|)$ implies that the action of $G/H$ is nonmodular. It then follows from Theorem~\ref{t1} that $\mathbb{K}[V]^{G/H}$ is a polynomial algebra. Consequently, the ring $\widehat{\mathbb{K}[W/H]^{G/H}_\mathfrak{m}} \simeq \widehat{\mathbb{K}[V]^{G/H}_\mathfrak{m}}$ is regular, as required.
 \hfill $\square$
\end{proof}

\section{Proof of the main theorem}
\label{s4}

In this section, we will prove Theorem~\ref{t4}.

Observe that a subgroup $H \subseteq G$ consisting of all unipotent matrices is the unique Sylow $p$-subgroup of the group $G$. In particular, $H$ is normal, therefore, proposition~\ref{p2} is applicable.

As was shown in the paper~\cite{10} of A. Braun, $p$-groups, generated by pseudoreflections in dimension 3, have polynomial ring of invariants:

\begin{theorem}[{\rm (see\  \cite{10})}]
\label{t6}
Let $H$ be a $p$-group, acting on a vector space $W$ of dimension 3 over a field $\mathbb{K}$ of characteristic $p$, generated by pseudoreflections. Then $\mathbb{K}[W]^H$ is polynomial.
\end{theorem}

Therefore, if $H$ is generated by pseudoreflections, then it follows from Theorem~\ref{t6} and Proposition~\ref{p2} that $\mathbb{K}[W]^G$ is polynomial (note that here $H$ may be non-abelian). Thus, from now on we will assume that $H$ is not generated by pseudoreflections. For abelian $H$ that are not generated by pseudoreflections we will show that $G \simeq A \times B$, where $A$ is a subgroup generated by homology-pseudoreflections of $G$, and $B$ is a subgroup generated by transvection-pseudoreflections of $G$. Moreover, we will show that both $A$ and $B$ have polynomial ring of invariants. Nevertheless, we will show that if $|A| > 1, |B| > 1$, then $G$ has non-polynomial ring of invariants. But $|A| > 1$ since $H$ is not generated by pseudoreflections, hence $G$ has polynomial ring of invariants if and only if $|B| = 1$, i.e., when $G$ does not contain transvections (and, in particular, coincides with $A$).

To begin, let us understand how upper triangular pseudoreflections look like in dimension 3. A brief case analysis shows that they can be one of three types:
\begin{equation}
\label{eq17}
	\begin{pmatrix}
	\lambda & a & b \\
	0 & 1 & 0 \\
	0 & 0 & 1
	\end{pmatrix},
\quad
	\begin{pmatrix}
	1 & 0 & b \\
	0 & 1 & a \\
	0 & 0 & \lambda
	\end{pmatrix},
\quad
	\begin{pmatrix}
	1 & (\lambda-1)b & ab \\
	0 & \lambda & a \\
	0 & 0 & 1
	\end{pmatrix},
\qquad \lambda \in \mathbb{K}^*, 
\quad a, b \in \mathbb{K}.
\end{equation}

Matrices of the first type will be called \textit{horizontal}, those of the second type \textit{vertical}, and those of the third type \textit{mixed}. If a pseudoreflection is both vertical and horizontal, it will be called \textit{corner}. Observe that a non-identity pseudoreflection, acting by a unitriangular matrix, is a transvection.

Since $H$ is not generated by pseudoreflections, it contains an element $g$ that is not a pseudoreflection. Then by a suitable choice of basis, we can ensure that the group $G$ also acts by upper triangular matrices and the distinguished element $g$ acts as the matrix
\begin{equation}
\label{eq18}
\begin{pmatrix}
1 & 1 & 0 \\
0 & 1 & 1 \\
0 & 0 & 1
\end{pmatrix}.
\end{equation}

Since $H$ is abelian by assumption, $H$ lies in the centralizer of the matrix $g$ in the group of all unipotent upper triangular matrices, which is
\begin{equation}
\label{eq19}
C_{U(W)}(g)=\left\{\begin{pmatrix}
1 & a & b \\
0 & 1 & a \\
0 & 0 & 1
\end{pmatrix}\Biggm|a, b \in \mathbb{K}\right\}.
\end{equation}

Next, by conjugating $g$ with a general pseudoreflection and considering that the result must lie in $H \subseteq C_{U(W)}(g)$, it is easy to see that it must be either a mixed pseudoreflection with $-1$ in the middle diagonal entry or a corner pseudoreflection. Moreover, multiplying two mixed pseudoreflections $\sigma, \sigma'$ with $-1$ in the middle diagonal entry yields an element of $H$, and from the abelian property it follows that its superdiagonal entries are equal. This is equivalent to the sums of the superdiagonal entries of $\sigma$ and $\sigma'$ being equal.

Thus, the group $G$ is generated by corner pseudoreflections and mixed pseudoreflections with $-1$ on the diagonal and a fixed sum $s$ of superdiagonal entries. Denote by $A$ the subgroup generated by mixed pseudoreflections, and by $B$ the subgroup generated by corner pseudoreflections. Let us examine the groups $A$ and $B$ in more detail and study their rings of invariants. As will be seen later, $A \cap B = \{1\}$, and also $B$ lies in the center of $G$, whence we have $G \simeq A \times B$.

First, let us study the group $A$. Let $\sigma_1, \ldots, \sigma_{t}$ be all the mixed pseudoreflections of the group $A$,
\begin{equation}
\label{eq20}
 \sigma_i=\begin{pmatrix}
 1 & s-c_i & \dfrac{c_i(c_i-s)}{2} \\
 0 & -1 & c_i \\
 0 & 0 & 1
 \end{pmatrix}.
\end{equation}
Then
\begin{equation}
\label{eq21}
 \sigma_i \sigma_j=\begin{pmatrix}
 1 & c_i-c_j & \dfrac{(c_i-c_j)(c_i-c_j-s)}{2} \\
 0 & 1 & c_i-c_j \\
 0 & 0 & 1
 \end{pmatrix} \in H.
\end{equation}

By direct computation, one can see that the group $A$ consists of the matrices $\sigma_1, \ldots, \sigma_t$ and the set of matrices
\begin{equation}
\label{eq22}
A_K=\left\{\begin{pmatrix}
1 & \alpha & \dfrac{\alpha(\alpha-s)}{2} \\
0 & 1 & \alpha \\
0 & 0 & 1
\end{pmatrix}\Biggm| \alpha \in K\right\},
\end{equation}
where $K$ is a subgroup of the additive group of the field $(\mathbb{K}, +)$, composed of the differences $c_i-c_j$, $i, j = 1, \ldots, t$. In particular, any pseudoreflection $\sigma_i$ can be represented as $\sigma_1 h$ or $h' \sigma_1$, where $h, h' \in A_K$. It follows that the group $A$ can be defined using the group $A_K$ and a single pseudoreflection $\sigma_1$. In particular, $|A| = 2p^k$, since $|A_K| = p^k$ for some non-negative integer $k$. Observe that if $k= \nobreak 0$, then the maximal unipotent subgroup $H$ of the group $G$ coincides with $B$ which is generated by pseudoreflections, whence $\mathbb{K}[W]^G$ is polynomial by Theorem~\ref{t6} and Proposition~\ref{p2}. Therefore, we will assume that $k > 0$.

Let us show that the ring of invariants of the group $A$ is polynomial. To do this, we define the norm of a polynomial under the group action:

\begin{definition}
\label{d4}
The effective norm $\operatorname{EN}_G(f)$ of a polynomial $f$ under the action of a finite group $G$ is the product
 \begin{equation}
 \label{eq23}
 \operatorname{EN}_G(f)=\prod_{g\operatorname{St}_G(f) \in G/\operatorname{St}_G(f)} g \cdot f,
 \end{equation}
where the product is taken over each coset in $G/\operatorname{St}_G(f)$, and an arbitrary representative is chosen from each class. From the definition, it is clear that $\operatorname{EN}_G(f) \in \mathbb{K}[W]^G$.
\end{definition}

The invariants of the group $A$ will be the polynomials
\begin{equation}
\label{eq24}
 f_1=z,
 \qquad f_2=y(cz-y)+2x'z,
 \qquad f_3=\operatorname{EN}_A(x'),
\end{equation}
where $c=c_1$, $x'=x+(s-c)y/2$.

The fact that $f_2 \in \mathbb{K}[W]^A$ can be verified by direct computation, and $f_3 \in \mathbb{K}[W]^A$ by the definition of the norm. Moreover, $|\operatorname{St}_A(x')|=2$, so $\deg f_3=\nobreak p^k$. The set of common zeros of the polynomials $f_1, f_2, f_3$ is $\{0\}$, and the product of their degrees equals the order of the group $A$, $|A| = 1 \cdot 2 \cdot p^k$. Therefore, by Lemmas~\ref{l2} and~\ref{l3}, the polynomials $f_1, f_2, f_3$ form a basis for the algebra of invariants of the group~$A$.

Since $x' \in \mathbb{K}[W]^{\sigma_1}$, $\operatorname{EN}_A(x')$ can be described as
\begin{equation}
\label{eq25}
 \operatorname{EN}_A(x')=\prod_{h \in A_K} h \cdot x'=\prod_{\alpha \in K} \biggl(x+\frac{s-c}{2}y+\alpha y+\frac{\alpha(\alpha-s)}{2}z+\dfrac{\alpha(s-c)}{2} z\biggr).
\end{equation}

In subsequent arguments, we will repeatedly substitute $z = 0$ into~$\operatorname{EN}_A(x')$. Note that in this case,
\begin{equation}
\label{eq26}
 \operatorname{EN}_A(x')|_{z=0}=\prod_{\alpha \in K} (x'+\alpha y)=\operatorname{EN}_{A'}(x'),
\end{equation}
where
\begin{equation}
\label{eq27}
 A'=\left\{
 \begin{pmatrix}
 1 & \alpha & 0 \\
 0 & 1 & 0 \\
 0 & 0 & 1
 \end{pmatrix}\Biggm| \alpha \in K \subseteq (\mathbb{K},+)\right\}
\end{equation}
is a certain group of horizontal matrices acting identically on $y$ and $z$. In what follows we will use the following technical proposition.

\begin{proposition}
\label{p3}
In the notations above the norm $\operatorname{EN}_{A'}(x')$ can be described as
\begin{equation}
\label{eq28}
	\operatorname{EN}_{A'}(x')=\sum_{i=0}^k \lambda_i y^{p^k-p^i}(x')^{p^i},
\qquad \lambda_i \in \mathbb{K},
\end{equation}
where $|A'|=p^k$ and $\lambda_0 \ne 0$.
\end{proposition}

\begin{proof}
Observe that the group $A'$ is isomorphic to a subgroup $K$ of the additive group of the field $\mathbb{K}$. Hence, $A'$ is a $k$-dimensional vector space over the prime subfield $\mathbb{Z}/p\mathbb{Z}$. Choose a basis $e_1, \ldots, e_k$ of $A'$ over $\mathbb{Z}/p\mathbb{Z}$. We prove the assertion by induction on $k$.

The base case $k = 0$ is immediate. Suppose that the assertion holds for $A''=\langle e_1, \dots, e_{k-1}\rangle$. We show that it also holds for $A'=\langle A'', e_k \rangle$.

Indeed, let
\begin{equation}
\label{eq29}
q(x',y)=\operatorname{EN}_{A''}(x')=\sum_{i=0}^{k-1} \mu_i y^{p^{k-1}-p^i}(x')^{p^i},
\qquad \mu_i \in \mathbb{K},
\quad \mu_0 \ne 0.
\end{equation}
Suppose that $e_k \cdot x'=x'+\alpha_k y$. Then, as is readily seen,
\begin{align}
\notag
e_k \cdot \operatorname{EN}_{A''}(x')
&=\sum_{i=0}^{k-1} \mu_i y^{p^{k-1}-p^i} (e_k \cdot (x')^{p^i})
=\sum_{i=0}^{k-1} \mu_i y^{p^{k-1}-p^i} (x'+\alpha_ky)^{p^i}
\\
&=\sum_{i=0}^{k-1} \mu_i y^{p^{k-1}-p^i} ((x')^{p^i}+\alpha_k^{p^i}y^{p^i})=q(x',y)+ q(\alpha_k y, y).
\label{eq30}
\end{align}

We claim that $q(\alpha_k y, y)=ry^{p^{k-1}} \ne 0, r \in \mathbb{K}$. By Lemmas~\ref{l2} and~\ref{l3}, the algebra $\mathbb{K}[W]^{A''}$ is a polynomial algebra generated by $q(x', y), y, z$. Suppose that $ry^{p^{k-1}}=0$. Then $q(x', y)$ is invariant under the action of $e_k$, and therefore $q(x', y) \in \mathbb{K}[W]^{A'}$. It follows that $\mathbb{K}[W]^{A''} \subset \mathbb{K}[W]^{A'}$. On the other hand, since $A''$~is a subgroup of $A'$, we have $\mathbb{K}[W]^{A'} \subset \mathbb{K}[W]^{A''}$. Thus, $\mathbb{K}[W]^{A'}= \mathbb{K}[W]^{A''}$. Consequently, $\mathbb{K}[W]^{A'}$~is a polynomial algebra generated by $q(x', y), y, z$. This contradicts Lemma~\ref{l3}: the polynomials $q(x', y), y, z$ form a homogeneous system of parameters and generate the algebra $\mathbb{K}[W]^{A'}$, whereas the product of their degrees is not equal to~$|A'|$. Therefore, $ry^{p^{k-1}} \ne 0$.

We now compute $\operatorname{EN}_{A'}(x')$. We have
\begin{align}
\notag
\operatorname{EN}_{A'}(x')
&=\prod_{j=0}^{p-1} (e_k^j \cdot \operatorname{EN}_{A''}(x'))=\prod_{j=0}^{p-1}(q(x',y)+j \cdot ry^{p^{k-1}})
\\
\notag
&=q(x',y)^p - q(x',y)(ry^{p^{k-1}})^{p-1}
\\
\notag
&=-\mu_0 r^{p-1}y^{p^k-1}(x')+\sum_{i=1}^{k-1} (\mu_{i-1}^p-r^{p-1}\mu_i)y^{p^k-p^i} (x')^{p^i} + \mu_{k-1}^p (x')^{p^k}
\\
&=\sum_{i=0}^k \lambda_i y^{p^k-p^i}(x')^{p^i}.
\label{eq31}
\end{align}
In particular, $\lambda_0=-\mu_0r^{p-1} \ne 0$. Thus, $\operatorname{EN}_{A'}(x')$ has the required form. \hfill $\square$
\end{proof}

Similarly, we present the basic invariants for the group $B$. It is easy to see that $B$ is a finite-dimensional vector space over $\mathbb{Z}/p\mathbb{Z}$ of some dimension $l$, so $|B| = p^l$. The basic invariants in this case will be the polynomials
\begin{equation}
\label{eq32}
f_1'=z,
\qquad f_2'=y,
\qquad f_3'=\operatorname{EN}_B(x').
\end{equation}
As in the proof of Proposition~\ref{p3}, one can show that $\operatorname{EN}_B(x')$ contains a monomial that is linear in $x'$, and that $\operatorname{EN}_B(x')$ has degree $p^l$. The fact that these polynomials are algebraically independent and generate the algebra $\mathbb{K}[W]^B$ follows from Lemmas~\ref{l2} and~\ref{l3}. Similarly, one obtains that a set of basic invariants for the invariant ring of the group $\langle B, \sigma_1 \rangle$ is given by
\begin{equation}
\label{eq33}
f_1''=z,
\qquad f_2''=y(cz-y),
\qquad f_3''=\operatorname{EN}_B(x').
\end{equation}

Finally, let us proceed to the proof of the theorem. By assumption, $H$ is not generated by pseudoreflections, so $|A_K| > 1$. Also, $H$ contains transvections, so $|B| > 1$. Thus, $k, l > 0$.

Assume that $G$ has a polynomial ring of invariants, with a basis formed by homogeneous invariants $h_1, h_2, h_3$. It is standard to verify that all linear invariants are of the form $\lambda z, \lambda \in \mathbb{K}$; in particular, among the generators there must be such a linear invariant, and we may assume that $h_1 = z$.

Then $\deg h_2 \cdot \deg h_3=2p^{k+l}$, so we may assume that $\deg h_2 = p^u$, $\deg h_3 = 2p^{k+l-u}$ for some $1 \le u \le k+l$.

Next, by considering three cases, we show that $u = \max(k, l)$ and one of $h_2, h_3$ cannot simultaneously be a polynomial in the generators of the ring of invariants of group $A$ and in the generators of the ring of invariants of group $\langle B, \sigma_1 \rangle$.

\textsl{Case} 1: $k > l$.
We show that in this case $u \ge k$. Indeed, otherwise the polynomial $h_2$, being a polynomial in $\operatorname{EN}_A(x')$, $y(cz-y)+2x'z$, $z$, cannot depend nontrivially on $\operatorname{EN}_A(x')$, so it depends only on $y(cz-y)+2x'z$ and $z$. But $y(cz-y)+2x'z$ is a polynomial of even degree, and $\deg h_2 = p^u$ is odd, so $h_2$ must be divisible by $z$, which is impossible because $h_1 = z$.

On the other hand, we show that $u$ cannot be greater than $k$, i.e., $u = k$. Suppose $u \ge k+1$. Then $2p^{k+l-u} \le 2p^{l-1} < p^l$, i.e., $h_3$ as a polynomial in $\operatorname{EN}_B(x')$,
$y(cz-y)$, $z$ cannot depend nontrivially on $\operatorname{EN}_B(x')$, whence it follows that $h_3$ is a polynomial in $y(cz-y)$ and $z$. On the other hand, it must be a polynomial in $\operatorname{EN}_A(x')$, $y(cz-y)+2x'z$, $z$, and since $\deg h_3 < \deg \operatorname{EN}_A(x')$, $h_3$ is a polynomial in $y(cz-y)+2x'z$ and $z$. But $h_3$ does not depend on $x'$, so $h_3$ can only depend on $z$, which contradicts algebraic independence.

In the case $u = k$, we find that $h_2$ can be expressed as a polynomial in the generators of $\mathbb{K}[W]^A$, and also as a polynomial in the generators of $\mathbb{K}[W]^{\langle B, \sigma_1\rangle}$, leading to the equality:
 \begin{equation}
 \label{eq34}
 \alpha \operatorname{EN}_A(x')+P\bigl(y(cz-y)+2x'z, z\bigr)
 =\alpha \operatorname{EN}_B(x')^{p^{k-l}}+Q\bigl(\operatorname{EN}_B(x'), y(cz-y), z\bigr),
 \end{equation}
where $Q$ does not contain $\operatorname{EN}_B(x')$ to the power $p^{k-l}$. The coefficients $\alpha$ in front of the norms on the left and right are equal because they are the coefficients of $(x')^{p^k}$ when considering $h_2$ as a polynomial in $x', y, z$.

We also see that the polynomial $P$ has total weighted degree $p^k$, and from the even degree of $y(cz-y)+2x'z$ it follows that $P$ is divisible by $z$. In particular, this implies that $\alpha \ne 0$, otherwise $h_2 = P$ would be divisible by $z$. Then we may assume $\alpha = 1$. Substituting $z = 0$ into the above equation and using formula \eqref{eq26}, we obtain
 \begin{equation}
 \label{eq35}
 \operatorname{EN}_{A'}(x')=(x')^{p^{k}}+Q_0((x')^{p^l}, -y^2),
 \end{equation}
    where $Q_0((x')^{p^l}, -y^2)=Q((x')^{p^l}, -y^2, 0)$.
But this equality is impossible because $\operatorname{EN}_{A'}(x')$ contains a monomial with $x'$ to the first power, while in $Q_0$ the variable $x'$ appears in all monomials with degree at least $p^l$. This concludes the analysis of case $k > l$.

 \textsl{Case} 2: $k=l$.

    Similarly to the previous case, it can be shown that $u = k = l$.

If $u = k$, then we have the equality
\begin{equation}
 \label{eq36}
\alpha \operatorname{EN}_A(x')+P(y(cz-y)+2x'z, z)=\alpha \operatorname{EN}_B(x')+Q(y(cz-y), z).
\end{equation}
Both $P$ and $Q$ are divisible by $z$, so $\alpha \ne 0$ and we may assume $\alpha = 1$. Substituting $z = 0$ into the above equation, we get
\begin{equation}
 \label{eq37}
\operatorname{EN}_{A'}(x')=(x')^{p^k},
\end{equation}
since $\operatorname{EN}_B(x')$ depends nontrivially only on $x'$ and $z$. But this equation is impossible because $\operatorname{EN}_{A'}(x')$ depends nontrivially on $y$.

\textsl{Case} 3: $k < l$.

In this case, similarly to the previous ones, it can be shown that $u = l$.

    In the case $u = l$, we have an equation similar to the first case, expressing $h_2$:
 \begin{equation}
  \label{eq38}
 \operatorname{EN}_A(x')^{p^{l-k}}+P(\operatorname{EN}_A(x'), y(cz-y)+2x'z, z)=\operatorname{EN}_B(x')+Q(y(cz-y), z).
 \end{equation}
    This time, $Q$ is divisible by $z$, and $P$ does not contain $\operatorname{EN}_A(x')$ to the power $p^{l-k}$. 
Substituting $z=0$ into this equation, we obtain
 \begin{equation}
  \label{eq39}
 \operatorname{EN}_{A'}(x')^{p^{l-k}}+P_0(\operatorname{EN}_{A'}(x'), -y^2)=(x')^{p^l},
 \end{equation}
    where $P_0(\operatorname{EN}_{A'}(x'), -y^2)=P(\operatorname{EN}_{A'}(x'), -y^2, 0)$.
    Here, the minimal degree of $x'$ among the monomials of $\operatorname{EN}_{A'}(x')^{p^{l-k}}$ is $p^{l-k}$, but since $P_0$ does not contain the term $\operatorname{EN}_{A'}(x')^{p^{l-k}}$, the minimal degree of $x'$ among the monomials of $P_0$ is less than $p^{l-k}$, so it cannot cancel out when added to $\operatorname{EN}_{A'}(x')^{p^{l-k}}$, making it impossible for the sum to equal $(x')^{p^l}$.

    The resulting contradiction concludes the analysis of the cases and the proof of the theorem.

\begin{remark}
\label{r3}
 In the above proof we used only that $H$ is abelian and $H$ contains an element that is not a pseudoreflection.
\end{remark}

The proven theorem allows us to construct a series of counterexamples in dimension 3 to the converse of Theorem~\ref{t1} in the reducible indecomposable case. Here is one of them:
    
\begin{example}
\label{e1}
	Upper triangular group $G$, generated by matrices
	\begin{equation}
 \label{eq40}
	\sigma_1=\begin{pmatrix}
	1 & 1 & 0 \\
	0 & -1 & 0 \\
	0 & 0 & 1
	\end{pmatrix},	
\qquad
	\sigma_2=\begin{pmatrix}
	1 & 0 & 0 \\
	0 & -1 & 1 \\
	0 & 0 & 1
	\end{pmatrix},
\qquad
	\tau=\begin{pmatrix}
	1 & 0 & 1 \\
	0 & 1 & 0 \\
	0 & 0 & 1
	\end{pmatrix},
	\end{equation}
	acting on a vector space of dimension 3 over a field of odd characteristic $p$, has a non-polynomial ring of invariants. Here, $G$ acts by a reducible indecomposable representation and $G \simeq D_p \times (\mathbb{Z}/p\mathbb{Z}),$ where $D_p$ is the dihedral group of order $2p$.
\end{example}

\begin{proof}
It is clear that $G$, like any upper triangular group, acts reducibly. On the other hand, if a subspace spanned by a vector $v$ is $G$-invariant, then the system of equations
$$
\begin{cases}
\sigma_1 \cdot v=\lambda_1 v, 
\\
\sigma_2 \cdot v=\lambda_2 v, 
\\
\tau \cdot v=\lambda_3 v 
\end{cases}
$$
readily implies that $\langle v \rangle = \langle e_1 \rangle$, where $e_1, e_2, e_3$ is a standard basis. Similarly, one can check that the only two-dimensional invariant subspace of $G$ is $\langle e_1, e_2 \rangle$. Therefore, the representation is indecomposable.

In the notation of the proof of Theorem~\ref{t4} we have $A=\langle \sigma_1, \sigma_2 \rangle, B=\langle \tau \rangle$. Moreover,
\begin{equation}
 \label{eq41}
\sigma_1 \sigma_2=\begin{pmatrix}
1 & 1 & 0 \\
0 & -1 & 0 \\
0 & 0 & 1
\end{pmatrix}
\begin{pmatrix}
1 & 0 & 0 \\
0 & -1 & 1 \\
0 & 0 & 1
\end{pmatrix}=
\begin{pmatrix}
1 & -1 & 1 \\
0 & 1 & -1 \\
0 & 0 & 1
\end{pmatrix}=r,
\end{equation}
whence the maximal unipotent subgroup $H=\langle r, \tau \rangle$ is abelian and contains transvections, but is not generated by them. It follows from Theorem~\ref{t4} that the algebra of invariants $\mathbb{K}[W]^G$ is not polynomial. The proof of Theorem~\ref{t4} also shows that $G \simeq A \times B$. Furthermore, $r^p=1$, $\sigma_1^2=1$ and $\sigma_1 r \sigma_1^{-1}= r^{-1}$. Thus, $A \simeq D_p$, while, as is readily seen, $B \simeq \mathbb{Z}/p\mathbb{Z}$. Consequently, $G \simeq D_p \times (\mathbb{Z}/p\mathbb{Z})$. \hfill $\square$
\end{proof}
\end{fulltext}

\end{document}